# COMPARISON OF THE METHODS FOR DISCRETE APPROXIMATION OF THE FRACTIONAL–ORDER OPERATOR

**Ľubomír DORČÁK[1], Ivo PETRÁŠ[2], Ján TERPÁK[1] and Martin ZBOROVJAN[1]**

[1] Department of Informatics and Process Control, BERG Faculty,
Technical University of Košice, B. Němcovej 3, 042 00 Košice, Slovak Republic
{*lubomir.dorcak, jan.terpak, martin.zborovjan*}*@tuke.sk*
[2] Steve's Electronic Services, 37 – 31255 Upper MacLure Rd.
Abbotsford, B.C., V2T 5N4, Canada, e-mail: *petras@telus.net*

**Abstract:** In this paper we will present some alternative types of discretization methods (discrete approximation) for the fractional-order (FO) differentiator and their application to the FO dynamical system described by the FO differential equation (FDE). With analytical solution and numerical solution by power series expansion (PSE) method are compared two effective methods - the Muir expansion of the Tustin operator and continued fraction expansion method (CFE) with the Tustin operator and the Al-Alaoui operator. Except detailed mathematical description presented are also simulation results. From the Bode plots of the FO differentiator and FDE and from the solution in the time domain we can see, that the CFE is a more effective method according to the PSE method, but there are some restrictions for the choice of the time step. The Muir expansion is almost unusable.

## 1 Introduction

The fractional-order calculus (FOCA) is about 300-years old topic. The theory of FO derivative was developed mainly in the 19th century. In the last decades besides theoretical research of the FO derivatives and integrals [Oldham, Spanier 1974, Samko, Kilbas, Marichev 1987, Podlubný 1994, 1999, 2000 and many others] there are growing number of applications of the FOCA in such different areas as e.g. long electrical lines, electrochemical processes, dielectric polarization [Westerlund 1994], colored noise, viscoelastic materials, chaos, in control theory [Manabe 1961, Oustalup 1988, 1995, Axtel, Bise 1990, Dorčák 1994, Podlubný, Dorčák, Koštial 1997, etc.] and in many other areas.

In works [Manabe 1961, Oustalup 1988, Axtel and Bise 1990, etc.] the first generalizations of analysis methods for FO control systems were made (s-plane, frequency response, etc.). Some of our works were oriented to the methods of FO system parameters identification, methods of FO controllers synthesis [Dorčák 1994, Podlubny et al. 1999, Petráš 1998, etc.], methods of stability analysis [Dorčák et al. 1998, Petráš et al. 1999], etc. The key step in digital implementation of the FO differentiator in the $PI^\lambda D^\delta$ controllers is their discretization.

## 2 Methods of discretization of the FO differentiator or integrator

In general, there are two discretization methods: direct discretization and indirect discretization. In indirect methods, two steps are required, i.e., frequency domain fitting in continuous time domain first and then discretizing the fit s-transfer function [Oustaloup et al. 2000, Vinagre et al. 2000]. In this paper, we will focus on the direct discretization methods. There are several approaches to direct discretization. We will consider only few of them. First is a power series expansion (PSE), second is continued fraction expansion (CFE), third is matrix approach proposed by [Podlubny 2000], fourth is a Taylor/MacLaurin series expansion of generating function suggested by [Duarte et al. 2002, Machado et al. 1997] and last one is Diethelm's approach based on solving of FDE in discrete time steps. All approaches are based on evolution of a generating function. As generating function we can use a simple Euler's backward difference, the Tustin rule, the Simpson rule or the Al-Alaoui operator, which is mixed scheme of Euler's rule and the Tustin trapezoidal rule. The investigation of these methods and comparison with each other and also with other methods was studied by several authors. The simplest and most straightforward method is the direct discretization using finite memory length expansion of Euler's rule from the Grunwald-Letnikov definition of fractional derivative. This approach was described by [Samko et al. 1987, Dorčák 1994, 2002, Gorenflo 1996, Podlubny 1999, etc.]. Gorenflo also suggested a PSE of the Tustin rule. The high order approximation by using the quadrature formula was used in [Lubich 1986] work. All these ways of approximation lead to form a FIR filter. It is well known that, for interpolation or evaluation purposes, rational functions are sometimes superior to polynomials, roughly speaking, because of their ability to model functions with zeros and poles [Vinagre et al., 2000, 2001] and therefore form of an IIR filter is much better. On the other hand we also do not need a large number of coefficients. This form of approximation was described by Vinagre and Chen and compared with PSE and the others as for example Muir recursion [Chen et al. 2002] or continuous methods [Vinagre et al. 2000]. Last but not least we should mention the approach proposed by Hwang [Hwang et al. 2002], which is based on B-splines function.

Muir recursion applied to the Tustin generating function (GF) and CFE applied to the Tustin or Al-Alaoui GF leads to the following approximations of the FO differentiator or integrator :

$$D^{\pm\delta}(z) = \left(\frac{K_1}{K_2 T}\right)^{\pm\delta} \lim_{n\to\infty}\left(\frac{A_n(z^{-1},+\delta)}{A_n(z^{-1},-\delta)}\right) = \left(\frac{K_1}{K_2 T}\right)^{\pm\delta} \frac{P_n(z^{-1},+\delta)}{Q_n(z^{-1},-\delta)} , \qquad (1)$$

$$D^{\pm\delta}(z) = \left(\frac{K_1}{K_2 T}\right)^{\pm\delta} \text{CFE}\left\{\left(\frac{1-z^{-1}}{1+z^{-1}/K_2}\right)^{\pm\delta}\right\}_{p,q} = \left(\frac{K_1}{K_2 T}\right)^{\pm\delta} \frac{P_p(z^{-1})}{Q_q(z^{-1})} , \qquad (2)$$

where $T$ is the sample period, $\delta$ is the order of differentiator or integrator, $K_1$ and $K_2$ are the constants according to the GF ($K_1=2$, $K_2=1$ for the Tustin GF, $K_1=8$, $K_2=7$ for the Al-Alaoui GF), $P$ and $Q$ are polynomials of degrees $p$ and $q$, respectively, in the variable $z^{-1}$ ($n=p=q$ for Muir recursion). For the numerical computation we have used $n=p=q=5$, 7 and 9. For e.g. $n=p=q=5$ the polynomials have the following coefficients for the Muir recursion : $Q_5(z^{-1})$ = $1/5\delta z^{-5}+1/5\delta^2 z^{-4}+(1/3\delta+1/15\delta^3)z^{-3}+2/5\delta^2 z^{-2}+\delta z^{-1}+1$, for CFE with the Tustin GF : $Q_5(z^{-1}) = (\delta^5-20\delta^3+64\delta)z^{-5}+(-195\delta^2+15\delta^4+225)z^{-4}+(105\delta^3-735\delta)z^{-3}+(420\delta^3-1050)z^{-2}-945\delta z^{-1}+945$, (the polynomials $P_9$ has the same coefficients except the odd coefficients, which have opposite sign) and for CFE with Al-Alaoui GF : $Q_5(z^{-1})=\{(1024\delta^5-11520\delta^4+40000\delta^3-31680\delta^2-51644\delta+51165)z^{-5}$ + $(26880\delta^4-282240\delta^3+920640\delta^2-882000\delta-92925)z^{-4}$ + $(329280\delta^3-2963520\delta^2+7696920\delta-5093550)z^{-3}+(2304960\delta^2-15558480\delta+23409750)z^{-2}+(9075780\delta-34034175)z^{-1}$+ 15882615} and $P_5(z^{-1})=\{(-1024\delta^5-11520\delta^4-40000\delta^3-31680\delta^2+51644\delta+51165)z^{-5}+(26880\delta^4+282240\delta^3+920640\delta^2+882000\delta-92925)z^{-4}$ + $(-329280\delta^3-2963520\delta^2-7696920\delta-5093550)z^{-3}+(2304960\delta^2+155584 80\delta+23409750)z^{-2}+(-90757 80\delta-34034175)z^{-1}+15882615$.

## 3 Comparison of the discretization methods

For all three above mentioned discretization methods we have compared their Bode characteristics first only for a simple FO differentiator and second by their application to the following FDE :

$$a_1 y^{(\delta)}(t) + a_0 y(t) = u(t) , \qquad (3)$$

where $\delta$ is order of differentiation - generally a real number, $a_1$, $a_0$ are arbitrary constants. This FDE (3) we have used only for evaluation purposes, because the FO $PI^\lambda D^\delta$ controller, which is described by the FDE $u(t) = Ke(t) + T_i D_t^{-\lambda} e(t) + T_d D_t^\delta e(t)$, after their application leads to the similar FDE - only more complex. For equation (3) we have compared the Bode characteristic too for their corresponding approximated transfer function :

$$G(z) = \frac{\sum_{i=0}^{5} Q_i\, z^{-i}}{a_0 \sum_{i=0}^{5} Q_i\, z^{-i} + a_1 \left(\frac{K_1}{K_2 T}\right)^{\pm\delta} \sum_{i=0}^{5} B_i\, z^{-i}} \qquad (4)$$

and we have made comparison of this approximation methods in the time domain too. We have compared the analytical solution of the FDE (3) [Podlubný 1994, Dorčák 1994] for $\delta \in (0,2\rangle$:

$$y(t) = \frac{1}{a_1} t^\delta E_{\delta,\delta+1}\left(-\frac{a_0}{a_1} t^\delta\right), \quad E_{\alpha,\beta}(z) = \sum_{k=0}^{\infty} \frac{z^k}{\Gamma(\alpha k + \beta)} \qquad (5)$$

(and for spatial cases $\delta=1$ or 2 the following classical analytical solutions of equation (3) too $y(t)=1/a_0(1-\exp(-a_0/a_1 t))$, $y(t)=1/a_0(1-\cos((a_0/a_1)^{1/2} t))$ with the numerical solution of the FDE (3) based on the Grunwald-Letnikov definition of the FO operator - PSE [Dorčák 1994, etc.]:

$$y_k = \frac{u_k - a_1 T^{-\delta} \sum_{i=1}^{k} b_j y_{k-j}}{a_1 T^{-\delta} + a_0}, \quad k = 1,2,\ldots \text{for } \delta = (0,1\rangle, \quad k = 2,3,\ldots \text{for } \delta = (1,2\rangle, \quad (6)$$

where $b_j$ are binomial coefficients $b_0^{(\delta)} = 1, b_j^{(\delta)} = (1 - (1+(\pm\delta))/j) b_{j-1}^{(\delta)}$, and with the numerical solution of the FDE (3) based on the Muir recursion or CFE with Tustin or Al-Alaoui operator :

$$y_k = \frac{1}{a_1 \left(\frac{K_1}{K_2 T}\right)^{\delta} P_0 + a_0 Q_0} \left[ -a_1 \left(\frac{K_1}{K_2 T}\right)^{\delta} \sum_{i=1}^{5} P_i y_{k-i} - a_0 \sum_{i=1}^{5} Q_i y_{k-i} + \sum_{i=0}^{5} Q_i u_{k-i} \right]. \quad (7)$$

We can derive analytical relations for the Bode characteristic of the FO differentiator $F_D(s) = T_D s^{\delta}$ :

$$\ln(F_{D^{\delta}}(i\omega)) = \ln(T_D) + \delta \ln(\omega) + i\delta \frac{\pi}{2} \quad (8)$$

and for the FDE (1) ($F(s) = 1/(a_1 s^{\delta} - a_0)$) :

$$\ln(F(i\omega)) = 0{,}5 \ln(\text{Re}^2 + \text{Im}^2) - i \ \text{arctg}\left(\frac{\text{Im}}{\text{Re}}\right), \quad (9)$$

where :

$$\text{Re} = \frac{a_1 \omega^{\delta} \cos(\delta\pi/2) + a_0}{[a_1 \omega^{\delta} \cos(\delta\pi/2) + a_0]^2 + [a_1 \omega^{\delta} \sin(\delta\pi/2)]^2} ,$$

$$\text{Im} = \frac{a_1 \omega^{\delta} \sin(\delta\pi/2)}{[a_1 \omega^{\delta} \cos(\delta\pi/2) + a_0]^2 + [a_1 \omega^{\delta} \sin(\delta\pi/2)]^2} \quad (10)$$

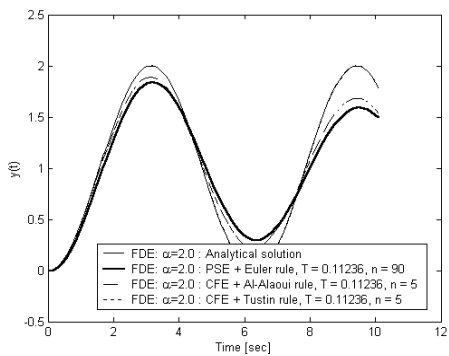 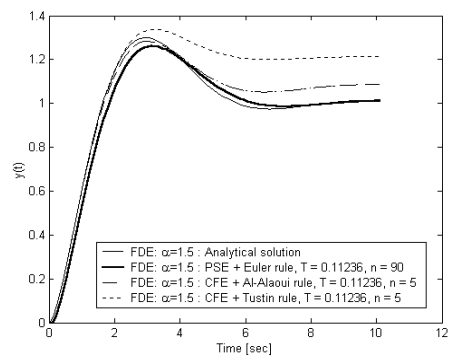

Figure 1.  Unit step responses       Figure 2.  Unit step responses

The numerical solution (PSE) (6) of the FDE (3) is in good agreement with the analytical solution (5) at $n > 90$ for all $\delta \in (0,2\rangle$. The agreement of this solution is very good for spatial case δ=1 or 2 with the classical analytical solutions. The accuracy was better by decreasing the time step. The numerical solution (7) with the Muir recursion was unusable for the time step near 0,1$s$, CFE with Tustin or Al-Alaoui GF are better (Fig. 1, 2). In these cases we can improve the accuracy by increasing of the $n$ to 7 or 9, but the accuracy was worse by decreasing of the time step. There is an border 0,1$s$ for the CFE and 0,5$s$ for the Muir recursion for all $n$=5, 7 or 9 for FDE (1). We can see from Fig. 3, 4, that the Bode plots of the discretized FO differentiator and FDE (1) are with large errors at higher frequencies and this discretization is unusable for smaller time step (unstable poles and/or zeros).

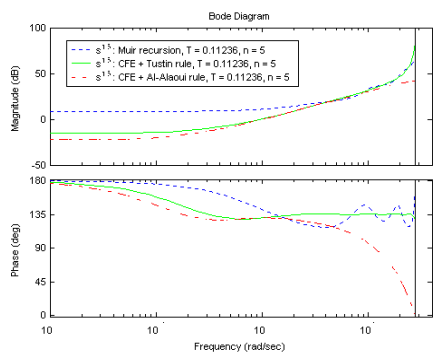 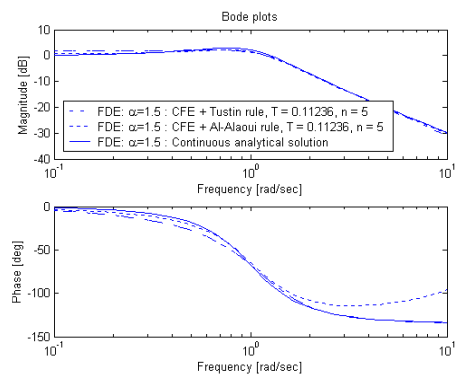

Figure 3. Bode plots of the differentiator      Figure 4. Bode plots of the FDE

This work was supported by grant VEGA 1/0374/03 from the Sl. Grant Agency for Science.